\documentclass[reqno,12pt]{amsart}

\usepackage[margin=2.5cm]{geometry}
\usepackage[T1]{fontenc}
\usepackage[utf8]{inputenc}
\usepackage[english]{babel}
\usepackage{amsmath,amssymb,amsthm,amsrefs}

\usepackage{enumitem}
\usepackage{mymacros}
\usepackage{hyperref}
\usepackage{xcolor}
\usepackage{orcidlink}
\usepackage{comment}
\usepackage{tikz}

\usepackage{centernot}
\usepackage{mathtools}

\DeclareMathOperator{\Hol}{Hol}
\DeclareMathOperator{\Har}{Har}
\DeclareMathOperator{\supp}{supp}

\usepackage{newpxtext,newpxmath}

\title{Random interpolation in the Nevanlinna and Smirnov classes and related spaces}

\author[Giuseppe Lamberti]{Giuseppe Lamberti \orcidlink{0009-0009-0503-0421}}
\address{Univ. Bordeaux, CNRS, Bordeaux INP, IMB, UMR 5251, F-33400 Talence, France} 
\email{giuseppe.lamberti@math.u-bordeaux.fr}

\begin{document}

\begin{abstract}
We study random interpolating sequences with prescribed radii in the Nevanlinna and Smirnov classes. As it turns out these are characterized by the Blaschke condition. This follows from a more general result. Indeed, we show that this characterization is true in so-called big Hardy-Orlicz spaces. It is noteworthy to mention that conditions for deterministic interpolation in these spaces are given by harmonic majorants, the existence of which is difficult to check in general.
\end{abstract}

\maketitle

\section{Introduction}
The Nevanlinna class is the most general class containing Hardy spaces and sharing many properties with these (Blaschke condition, factorization, non-tangential behavior, harmonic majorants, etc.), and interpolation in this class actually has some history. A first problem was stated and solved by Naftalevi\v{c} \cite{naftalevic56} in the late 50s, while Yanagihara \cite{yanagihara74} considered a related problem in the Smirnov class. 

In a more recent work \cite{hartmann2004}, the authors considered so-called free interpolation in the Nevanlinna and Smirnov classes and completely characterized the corresponding free interpolating sequences. The principal result they proved connects free interpolating sequences to harmonic majorants, which are notoriously difficult to determinate. The aim of this paper is to use a random model to derive a more geometric and "easy to check" condition for free interpolating sequences in the Nevanlinna and Smirnov classes, as has been done in other contexts \cites{chalmoukis22, bomash92}.

Before going into details, we recall some terminology that we will use.
Consider the unit disk $\bD = \{z \in \bC : |z|<1\}$ and its topological boundary $\bT := \partial \bD$. The \textit{Blaschke factor} at $\lambda \in \bD$ is
$$
\varphi_\lambda(z) := \frac{z-\lambda}{1-\ol{\lambda}z}, \quad w \in \bD,
$$
while the normalized Blaschke factor is
\begin{align*}
b_\lambda(z) & := \frac{|\lambda|}{\lambda} \frac{z-\lambda}{1-\ol{\lambda}z} \quad \lambda \neq 0,\\
b_0(z) & := z.
\end{align*}
Given $\lambda,\mu \in \bD$ define the \textit{pseudohyperbolic distance} between $\lambda$ and $\mu$ as 
$$
\rho(\lambda,\mu) := |\varphi_\lambda(\mu)|.
$$

The \textit{Nevanlinna class} is defined by
$$
\cN := \Big\{ f \in \Hol (\bD) : \lim_{r \to 1} \frac{1}{2\pi} \int_0^{2\pi} \log^+|f(re^{i\theta})| < \infty  \Big\}
$$
and the Smirnov class is given as:
$$
\cN^+ := \Big\{ f \in \cN : \lim_{r \to 1} \frac{1}{2\pi} \int_0^{2\pi} \log^+|f(re^{i\theta})| = \frac{1}{2\pi} \int_0^{2\pi} \log^+|f(e^{i\theta})| < \infty  \Big\}.
$$
We now describe the interpolation problems mentioned before. Let $X$ be a space of holomorphic functions on the disk and $\Lambda$ a sequence of elements contained in $\bD$. Consider the abstract trace space
$$
X|_\Lambda := \{ (f(\lambda))_{\lambda \in \Lambda}: f \in X \}.
$$
A sequence $\Lambda$ is \textit{free interpolating} for $X$ if
$$
\ell^\infty X|_\Lambda \subset X|_\Lambda. 
$$
It is not difficult to see that if $X$ is an algebra that contains constants (e.g. the Nevanlinna class) the condition is equivalent to
$$
\ell^\infty \subset X|_{\Lambda}.
$$
Free interpolating sequences have been widely studied in the last years, especially for their connection with unconditional bases in model spaces. For further details the interested reader can see \cite{Nikolski2002b}*{Part C}.
Observe that if $\Lambda$ is free interpolating for the Nevanlinna class, then it is a \textit{zero set}, namely it must exist a function $f\in \cN$ such that $f|_{\Lambda} = 0$ but $f$ is not identically zero. Note that zero sets for the Nevanlinna class are precisely \textit{Blaschke sequences}, i.e. sequences that satisfy
\begin{equation} \label{Blaschke sequence}
\sum_{\lambda \in \Lambda} (1-|\lambda|) < \infty
\end{equation}
(see for instance \cite{garnett06}). As a consequence, we obtain that every free interpolating sequence for $\cN$ satisfies \eqref{Blaschke sequence}. It is well known that whenever \eqref{Blaschke sequence} is satisfied we can define the corresponding Blaschke product:
$$
B_\Lambda(z) :=  \prod_{\lambda \in \Lambda} b_\lambda(z)
$$
where $\lambda$ appears according to its multiplicity. In the case of (simple) interpolation, all zeros of $B_\Lambda$ are simple, and, in this situation, we will also consider the Blaschke product with an element removed:
$$
B_\lambda(z) := B_{\Lambda \setminus \{\lambda\}} = \prod_{\mu \neq \lambda} b_\mu(z).
$$
In \cite{hartmann2004} the authors provided a complete characterization of free interpolating sequences for the Nevanlinna and Smirnov classes, using the function $\varphi_\Lambda : \bD \to \bR$ defined as
$$
\varphi_\Lambda (z) = 
\begin{dcases}
    \log \frac{1}{|B_\lambda(\lambda)|} & \text{if }z=\lambda \in \Lambda \\
    0 & \text{otherwise}.
\end{dcases}
$$
As it turns out, the characterization of deterministic interpolating sequences is related to harmonic majorants. Given a finite measure $\mu$ on $\bT$, we define the harmonic function
$$
P[\mu] (z) := \int_{\bT} \frac{1-|z|^2}{|\zeta-z|^2} d\mu(\zeta).
$$
A function $u$ defined on $\bD$ is said to admit a positive harmonic majorant if there is a positive finite measure $\mu$ on $\bT$ such that $u\leq P[\mu]$ on $\bD$. If $d\mu=w dm$, with $w \in L^1(\bT)$, we write $P[w]$ instead of $P[wdm]$ and we say that $u$ admits a positive \textit{quasi-bounded} harmonic majorant.
The next theorem is one of the main results proved in \cite{hartmann2004}.

\begin{thm} \label{free int Nevanlinna}
    Let $\Lambda$ be a Blaschke sequence. The following are equivalent
    \begin{enumerate}[label=(\alph*)]
        \item $\Lambda$ is interpolating for the Nevanlinna (resp. Smirnov) class.
        \item $\varphi_\Lambda$ admits a positive harmonic majorant (resp. quasi-bounded harmonic majorant).
        \item The trace space is given by:
        $$
        \cN\left|_{\Lambda} \right. = l_{\cN} := \left\{ (a_\lambda)_\lambda : \exists \, h \in \Har^+(\bD) \tst h(\lambda) \geq \log^+|a_\lambda|, \, \lambda \in \Lambda  \right\}.
        $$
        (respectively: 
        $$
        \cN^+\left|_{\Lambda} \right. = l_{\cN^+} := \left\{ (a_\lambda)_\lambda : \exists \, h \in \Har^+(\bD) \text{ quasi bounded} \tst h(\lambda) \geq \log^+|a_\lambda|, \, \lambda \in \Lambda  \right\}).
        $$
    \end{enumerate}
\end{thm}

On the other hand, Naftalevi\v{c} completely characterized the sequences $\Lambda$ such that
$$
\cN |_\Lambda = l_{Na} :=  \{(a_\lambda)_{\lambda \in \Lambda}: \sup_{\lambda \in \Lambda} (1-|\lambda|)\log^+|a_\lambda| < \infty\},
$$
proving that such a condition is reached if and only if $\Lambda$ is contained in a finite union of Stolz angles and
$$
\sup_{\lambda \in \Lambda} (1-|\lambda|)\log \frac{1}{|B_\lambda(\lambda)|} < \infty.
$$
As is well known, Hardy  space interpolating sequences do not need to be contained in a finite union of Stolz angles. In a sense, the trace space $l_{Na}$ does not seem to be an appropriate choice for $\cN$. In addition, we note that the random model we consider does almost surely not produce interpolating sequences in the sense of Naftalevi\v{c} (see Proposition \ref{naftalevic:random}).  For these reasons, we focus on free interpolating sequences.

We shall now present the random model we are interested in. If not differently specified, by the term \textit{random sequence} we mean a sequence of the following kind. Consider a family of fixed radii $(r_n)_{n \in \bN}$, with $0\leq r_n < 1$, $\sum (1-r_n) < \infty$ (Blaschke condition), and a family of random variables $(\theta_n)_{n \in \bN}$ all independent and uniformly distributed in $[0,2\pi)$ . Define 
$$
\Lambda(\omega) = (\lambda_n(\omega))_{n \in \bN} = (r_n e^{i \theta_n (\omega)})_{n \in \bN}.
$$
In this scenario we want to find necessary and sufficient conditions on the sequence $(r_n)_{n \in \bN}$ such that the corresponding random sequence is almost surely free interpolating for $\cN$ or $\cN^+$.

These random sequences, sometimes called Steinhaus sequences in the literature, have already been studied in the last 30 years. They were used by Cochran in \cite{cochran90}, where the author found a $0-1$ law for weakly separated sequences (see Section 2 for precise definitions), while Rudowicz \cite{rudowicz94} further contributed characterizing interpolating sequences for the Hardy space. We also recall \cites{bomash92, bogdan96} where this model was used to better understand zero sets in Dirichlet and Bergman spaces and \cite{massaneda96}, where the work of Cochran and Rudowicz was extended in the unit ball. It is worth mentioning more recent works, such as \cite{chalmoukis22}, which deals with interpolation in Dirichlet spaces, \cite{dayan2023}, where the authors considered Besov and Hardy spaces in several variables, both for the polydisc and the unit ball, \cite{chalmoukis2024} on random Carleson sequences on the polydisc and random interpolation in the Bloch space, and \cite{kononova2022} which deals with zero sets for Fock-type spaces on the complex plane.

The principal result of this paper is the following theorem.

\begin{thm} \label{THM}
    Let $\Lambda$ be a random sequence. Then $\Lambda$ is almost surely free interpolating for the Nevanlinna or the Smirnov classes if and only if it is a Blaschke sequence.
\end{thm}

As it turns out, this theorem follows from a much stronger result. Indeed, the harmonic majorant $P[\mu]$ of the function $\varphi_\Lambda$ introduced above and that we will construct later satisfies almost surely $d\mu=w dm$ with $w\in L^p(\bT)$, $p\geq 1$. Consequently, we do not only get almost sure interpolation in Nevanlinna and Smirnov classes but in smaller so-called  Hardy-Orlicz spaces as we will explain later.

The theorem can be reformulated saying that almost every random sequence which is a zero set for $\cN$ (or $\cN^+$) is a free interpolating sequence. It is worth pointing out that this is not the first time that such a scenario occurs. In \cite{chalmoukis22} a similar situation has been revealed regarding a range of weighted Dirichlet spaces.

\subsection{Outline of the main result} As mentioned above, in order to prove our main theorem, we will establish a stronger result. As stated in Theorem \ref{main: hardy-orlicz}, it turns out that almost every random sequence which is a Blaschke sequence is automatically free interpolating for a range of spaces that are contained in the Nevanlina and Smirnov classes, namely \textit{big Hardy-Orlicz spaces}, which we will introduce now. Let $\varphi: \bR \to [0,\infty) $ be a convex, non-decreasing function satisfying
\begin{enumerate}[label=(\roman*)]
    \item $\lim_{t \to \infty} \varphi(t)/t = \infty$
    \item  $\varphi(t+2)\leq M\varphi(t) + K$, $t\geq t_0$, for some constants $M,K \geq 0$ and $t_0 \in \bR$.
\end{enumerate}
Whenever $\varphi$ satisfies these two requirements, we say that $\varphi$ is strongly convex. If this is the case we can define the associated Hardy-Orlicz space as
$$
\cH_{\varphi \circ \log} := \Big\{ f \in \cN^+ : \int_\bT \varphi(\log|f(\zeta)|) dm(\zeta) < \infty \Big\}.
$$
It is worth mentioning that these spaces include Hardy spaces ($\varphi(t)=e^{pt}$), while the limit case $\varphi(t) = t$ (which is not a strongly convex function) corresponds to the Smirnov class $\cN^+$.

The study of free interpolating sequences for $\cH_{\varphi \circ \log}$ was conducted in \cite{hartmann07} which was subsequent to \cite{hartmann2004}. The spaces in question are Hardy-Orlicz spaces with a defining function $\varphi$ that must fulfill the following quasi-triangular inequality
\begin{equation} \label{quasi-triangular}
\varphi(a+b) \leq c(\varphi(a)+\varphi(b))
\end{equation}
for some constant $c\geq 1$ and for all $a,b \geq t_0$. This condition implies that $\cH_{\varphi \circ \log}$ is an algebra and does not correspond to any classical Hardy spaces $H^p$, $p>0$.
We report the main theorem \cite{hartmann07}*{Theorem 1.2}.
\begin{thm}
    Let $\varphi : \bR \to [0,\infty)$ be a strongly convex function satisfying \eqref{quasi-triangular}. The following assertions are equivalent
    \begin{enumerate} [label=(\alph*)]
        \item $\Lambda$ is free interpolating for $\cH_{\varphi \circ \log}$
        \item There exists a positive measurable function $w \in L^\varphi(\bT)$ such that $\varphi_\Lambda \leq P[w]$.
        \item The trace space is given by
        $$
        \cH_{\varphi \circ \log}|_\Lambda = l_{\varphi \circ \log} := \left\{ (a_\lambda)_\lambda : \exists \, w \in L^\varphi \tst \log^+|a_\lambda| \leq P[w], \, \lambda \in \Lambda  \right\}.
        $$
    \end{enumerate}
    where the space $L^\varphi(\bT)$ is the space of measurable functions $u$ such that $\varphi \circ |u| \in L^1(\bT)$.
\end{thm}

We now restrict ourselves to a specific range of spaces, namely the ones with defining function $\varphi(t)=t^p$ for $p > 1$. It is immediate that these functions are all strongly convex and satisfy the quasi-triangular inequality. For simplicity of notation we write $\cH_p$ instead of $\cH_{\varphi \circ \log}$. It is worth mentioning that in this scenario $L^\varphi = L^p$, $(L^\varphi)^*= L^q$ and furthermore that:  $\cH_q \subset \cH_p$, for $p<q$. Lastly we recall the following inclusions between the spaces of holomorphic functions we are considering:
$$
H^q \subset \cH_p \subset \cN^+ \subset \cN \quad \text{for every $p > 1, q>0$.}
$$
We point out that the Blaschke condition characterizes zero sets in $H^q, \cH_p$, $\cN$ and $\cN^+$.

We are now ready to state the theorem regarding the spaces $\cH_p$.

\begin{thm} \label{main: hardy-orlicz}
    Let $\Lambda$ be a random sequence. Then $\Lambda$ is almost surely free interpolating for $\cH_p$ if and only if it is a Blaschke sequence.
\end{thm}

Theorem \ref{THM} follows immediately from Theorem \ref{main: hardy-orlicz} as follows.

\begin{proof}[Proof of Theorem \ref{THM}]
Since every free interpolating sequence is a zero set and zero sets for the Nevanlinna or Smirnov classes are Blaschke sequences, the necessary part follows.

    Now let $\Lambda$ be a random sequence which satisfies the Blaschke condition. From Theorem \ref{main: hardy-orlicz} we know that $\Lambda$ is almost surely free interpolating for $\cH_p$.
    Since $\cH_p \subset \cN^+ \subset \cN$ for every $p > 1$, this means that $\cH_p |_{\Lambda} \subset \cN^+|_\Lambda \subset \cN|_\Lambda$ for every $\Lambda$. Hence if $\Lambda$ is almost surely free interpolating for $\cH_p$ then
    $$
    \ell^\infty \subset \cH_p|_\Lambda \subset \cN^+|_\Lambda \subset \cN|_\Lambda \quad \text{almost surely,}
    $$
    which means that $\Lambda$ is almost surely free interpolating for the Nevanlinna and Smirnov classes.
\end{proof}

The paper is structured as follows. Section 2 contains the notions needed both on the complex analysis and probabilistic side. Section 3 is devoted to the proof of Theorem \ref{main: hardy-orlicz} and to investigate how Naftalevi\v{c}'s interpolation problem changes using the random model we considered.

\section{Preliminaries}

\subsection{Spaces of holomorphic functions}
We start this section recalling some properties that are of interest for our study. A sequence $\Lambda=(\lambda)_{\lambda \in \Lambda} \subset \bD$ is
\begin{itemize}[label=\labelitemii]
    \item \textit{weakly separated} if 
    $$
    \inf_{\mu \neq \lambda} \rho(\lambda,\mu) > 0;
    $$
    \item \textit{strongly separated} if
    $$
    \inf_{\lambda \in \Lambda} \prod_{\mu \neq \lambda} \rho(\lambda,\mu) > 0.
    $$
\end{itemize}

Let $\zeta \in \bT$ and $\alpha >1$, the \textit{Stolz angle} is
$$
\Gamma_\alpha(\zeta) := \{ z \in \bD: |\zeta-z| \leq \alpha (1-|z|)\}.
$$

Given a positive Borel measure $\mu$ on the disk, we say that $\mu$ is a \textit{Carleson measure} for the Hardy space $H^2$ if there exists $C=C(\mu)>0$ such that for every $f \in H^2$
$$
\int_{\bD} |f(z)|^2 d\mu(z) \leq C \|f\|^2_{H^2}.
$$
When $\mu=\mu_\Lambda := \sum_{\lambda \in \Lambda} (1-|\lambda|)\delta_\lambda$ is a Carleson measure, we say that $\Lambda$ is a \textit{Carleson sequence}. We recall that Carleson sequences in the sense we have just introduced are finite unions of $H^\infty$-interpolating sequences. As proved in \cite{carleson58} and \cite{carleson62}, all these notions are related to the interpolation problem Carleson solved for $H^\infty$.

We also introduce the $\textit{balayage}$ of a positive Borel measure $\mu$ on $\bD$:
$$
B\mu (\zeta) := \int_{\bD} \frac{1-|z|^2}{|\zeta-z|^2} d\mu(z).
$$

We mention the following useful result, contained in \cite{hartmann2004}*{Proposition 4.1}.

\begin{thm} \label{Nicolau_Prop}
    Let $\Lambda$ be a Blaschke sequence. For any $\delta \in (0,1)$ there exists a quasi-bounded positive harmonic function $h=P[w]$, $w \in L^1(\bT)$, such that
    $$
    \sum_{\mu: \rho(\mu,\lambda) \geq \delta} \log \frac{1}{|b_\lambda(\mu)|} \leq h(\lambda).
    $$
\end{thm} 
From this theorem it follows immediately that every weakly separated sequence is free interpolating for $\cN^+$ (and so for $\cN$).

On the Hardy-Orlicz side, we underline that the functions $\varphi(t)=t^p$ satisfy for all $p>1$ the so called $\nabla_2$-condition (see for instance \cite{lesniewicz73}). Specifically, there exist $c>1$ and $t_0 >0$ such that the following inequality holds:
$$
2\varphi(t) \leq \frac{1}{c} \varphi(ct), \quad t \geq t_0.
$$

With this in mind we can recall the following theorem, proved in \cite{hartmann07}*{Theorem 1.3}.

\begin{thm} \label{hardy orlicz thm 1.3}
    Let $\varphi$ be a strongly convex function that satisfies \eqref{quasi-triangular} and the $\nabla_2$-condition. If $u$ is a non-negative Borel function on the unit disk, the following two assertions are equivalent.
    \begin{enumerate}[label=(\alph*)]
        \item there exists a function $w \in L^\varphi$ such that $u(z) \leq P[w](z)$ for all $z \in \bD$.
        \item There exists $C\geq 0$ such that
        $$
        \sup_{\mu \in \cB_{\varphi^*}} \int_\bD u(z) d\mu(z) \leq C,
        $$
        where $\cB_{\varphi^*} = \{ \mu \text{ positive measure on } \bD: \|B\mu\|_{(L^\varphi)^*} \leq 1 \}$.
    \end{enumerate}
\end{thm}

\subsection{Probability background}

On the probabilistic side, we only need the following inequality, proved in \cite{rosenthal70}*{Lemma 1}.

\begin{thm}[Rosenthal's inequality] \label{rosenthal ineq}
    Let $1<p < \infty$ and consider $X_1,\ldots,X_k$ independent and non-negative random variables such that $\bE[X_m^p]<\infty$, for every $m=1,\ldots,k$. Then
    \begin{equation} \label{rosenthal}
        \bE\Big[ \big( \sum_{m=1}^k X_m \big)^p \Big] \leq 2^{p^2} \max\Big\{ \sum_{m=1}^{k} \bE[X_m^p],  \big( \sum_{m=1}^{k} \bE[X_m]\big)^p \Big\}.
    \end{equation}
\end{thm}
 In the literature there are many versions and generalizations of \eqref{rosenthal}, but for our purpose the version reported here is sufficient.
 
\subsection{Notation}
If $f$ and $g$ are positive expressions, we will write $f \lesssim g$ if there exists $C>0$ such that $f\leq C g$, where $C$ does not depend on the parameters behind $f$ and $g$. We will simply write $f \simeq g$ if $f \lesssim g$ and $g \lesssim f$.

\section{Proof of the main result}

In what follows we will strongly use the dyadic partition of the unit disk $\bD$. Consider for $n \geq 0$,
\begin{align*}
A_n & = \left\{ z \in \bD : 1-2^n \leq |z| < 1-2^{-(n+1)} \right\}, \\ 
N_n & = \# (A_n \cap \Lambda).
\end{align*}

Note that $N_n$ is actually a number and not a random variable, since it depends only on the sequence of the radii. 
We now reformulate the Blaschke condition using the notation just introduced. Note that
$$
\sum_{\lambda \in \Lambda} (1-|\lambda|) = \sum_{n \geq 0} \sum_{\lambda \in A_n} (1-|\lambda|) \simeq \sum_{n \geq 0} N_n 2^{-n},
$$
Hence the Blaschke condition is equivalent to
$$
\sum_{n \geq 0} N_n 2^{-n} < \infty.
$$
We emphasize that even though we are interested in random sequences, the property of being or not a Blaschke sequence (and therefore a zero set for $\cN$) does not involve probability. A major difference with the situation in, for instance, Dirichlet spaces, is that there is no useful deterministic description of zero sets known for the Dirichlet space; still, and somehow surprisingly, Bogdan \cite{bogdan96} was able to give an easy condition ensuring that a sequence is a.s. a zero set for the Dirichlet space.

The following result is standard. We produce here a proof for convenience of the reader.

\begin{lem} \label{poisson Lq norm}
    Let $P_z$ be the Poisson kernel on $\bD$ and $q\geq 1$. Then
    $$
    \|P_z\|_{L^q}^q \simeq \frac{1}{(1-|z|)^{q-1}}.
    $$
\end{lem}

\begin{proof}
    We have that
\begin{align*}
    \|P_z\|_{L^q}^q & = \int_\bT \frac{1-|z|^2}{|\zeta-z|^2} \left(\frac{1-|z|^2}{|\zeta-z|^2} \right)^{q-1} dm(\zeta) \\
    & \leq C \left(\frac{1}{1-|z|}\right)^{q-1}.
\end{align*}
On the other hand, if $z=|z|e^{i\theta_0}$ and $\zeta = e^{i\theta}$, then
$$
P_z(\zeta) = P_z(e^{i\theta})=\sum_{n \in \bZ} |z|^{|n|} e^{in(\theta-\theta_0)}.
$$
From Parseval's theorem it follows 
$$
\|P_z\|_{L^2}^2 = \sum_{n \in \bZ} |z|^{2|n|} = 2 \frac{1}{1-|z|^2},
$$
and by H\"older's inequality we have that
$$
\frac{2}{1-|z|^2} = \|P_z \cdot P_z\|_{L^1} \leq \|P_z\|_{L^q} \|P_z\|_{L^p} \leq \|P_z\|_{L^q} \frac{1}{(1-|z|)^{\frac{p-1}{p}}}
$$
which implies
$$
\|P_z\|_{L^q} \geq \frac{2}{1+|z|} \frac{1}{(1-|z|)^{\frac{q-1}{q}}},
$$
concluding the proof.
\end{proof}

\begin{lem} \label{alpha_lambda estimate}
    Let $\mu = \sum_{\lambda} \alpha_\lambda \delta_\lambda$ be a positive Borel measure on $\bD$ such that $\|B\mu\|_{L^q} \leq 1$, $q \geq 1$, then there exists $C>0$ such that
    
    \begin{equation} \label{alpha_lambda equation}
    \sum_{\lambda \in \Lambda} \alpha_\lambda^q \frac{1}{(1-|\lambda|)^{q-1}} \leq C.
    \end{equation}
    
\end{lem}

\begin{proof}
    Suppose $\|B\mu\|_{L^q} \leq 1$, then 
    \begin{align*}
        1 & \geq \int_\bT B\mu(\zeta)^q dm(\zeta) \\
        & \geq \sum_{\lambda \in \Lambda} \alpha_\lambda^q \int_\bT \left(\frac{1-|\lambda|^2}{|\zeta-\lambda|^2}\right)^q dm(\zeta).
    \end{align*}
    From Lemma \ref{poisson Lq norm}, we obtain the statement.
\end{proof}

\begin{prop} \label{generalization H_p}
    Let $\Lambda$ be a sequence in the unit disk and consider $p >1$. If
    $$
    \sum_{\lambda \in \Lambda} (1-|\lambda|) \log^p\frac{1}{|B_\lambda(\lambda)|} < \infty
    $$
    then $\Lambda$ is free interpolating for $\cH_p$.
\end{prop}

\begin{proof}
We want to apply Theorem \ref{hardy orlicz thm 1.3} to the function $\varphi_\Lambda$. Since $\varphi_\Lambda$ is different from zero only on $\Lambda$, it is enough to show that under our hypothesis
$$
\sup_{\substack{\mu : \|B\mu\|_{L^q} \leq 1 \\ \supp \mu \subset \Lambda}} \int_\bD \varphi_\Lambda(z) d\mu(z) < \infty. 
$$
Consider $\mu = \sum_\lambda \alpha_\lambda \delta_\lambda$. Using H\"{o}lder's inequality and Lemma \ref{alpha_lambda estimate} we have
\begin{align*}
    \int_\bD \varphi_\Lambda(z) d\mu(z) & = \sum_{\lambda \in \Lambda} \alpha_\lambda \varphi_\Lambda(\lambda) \\
    & \leq \left(\sum_{\lambda \in \Lambda} \alpha_\lambda^q \frac{1}{(1-|\lambda|)^{q-1}}\right)^{1/q} \left(\sum_{\lambda \in \Lambda} (1-|\lambda|)\log^p \frac{1}{|B_\lambda(\lambda)|}\right)^{1/p} \\
    & \leq C \left(\sum_{\lambda \in \Lambda} (1-|\lambda|)\log^p \frac{1}{|B_\lambda(\lambda)|}\right)^{1/p}.
\end{align*}
Since the last term is finite by hypothesis and does not depend on the choice of $(\alpha_\lambda)_{\lambda \in \Lambda}$, we obtain the result.
\end{proof}

\begin{rem}
Proposition \ref{generalization H_p} is an extension of \cite{hartmann2004}*{Corollary 1.11}, which deals with the case $p=1$. We report the proof for completeness. Consider the function
\begin{equation} \label{def psi}
\psi(\zeta) = \sum_{\lambda \in \Lambda} \log \frac{1}{|B_\lambda(\lambda)|} \chi_{I_\lambda}(\zeta).
\end{equation}
Suppose 
\begin{equation} \label{sum Blambda}
   \sum_{\lambda \in \Lambda} (1-|\lambda|) \log\frac{1}{|B_\lambda(\lambda)|} < \infty, 
\end{equation}
then $\psi \in L^1(\bT)$, so we can define its Poisson extension, $P[\psi]$, which is a harmonic majorant of $\varphi_\Lambda$. The result follows from Theorem \ref{free int Nevanlinna}.
\end{rem}

\begin{lem} \label{diagonal estimate}
    Let $\lambda=re^{i\theta_1}$, $\mu=se^{i\theta_2}$, $\theta_1,\theta_2$ independent and uniformly distributed in $[0,2\pi)$ with $1/2 \leq (1-r)/(1-s) \leq 2$ and $1/2 \leq r,s<1$. Then for every $p>0$ there exists a constant $C_p>0$ such that
    $$
    \bE\Big[\log^p\frac{1}{\rho(\lambda,\mu)}\Big] \leq C_p (1-r^2).
    $$
\end{lem}

\begin{proof}
Due to the fact that $\lambda$ and $\mu$ are independent and $1-r \simeq 1-s$, we have
\begin{align*}
\bE\Big[\log^p \frac{1}{\rho(\lambda,\mu)} \Big] & = \frac{2^{-p}}{4\pi^2} \int_0^{2\pi}\int_0^{2\pi} \log^p \frac{1}{\rho(re^{i\theta_1},se^{i\theta_2})^2} d\theta_1 d\theta_2 \\
& \lesssim \frac{2^{-p}}{2\pi} \int_0^{2\pi} \log^p \frac{1}{\rho(r,re^{i\theta})^2} d\theta
\end{align*}
Notice that
$$
\frac{1}{\rho(r,re^{i\theta})^2} = 1 + \frac{(1-r^2)^2}{2r^2(1-\cos\theta)} =  1 + \frac{(1-r^2)^2}{(r\theta)^2}(1+o(\theta)),
$$
which yields
\begin{align*}
\bE\Big[\log^p \frac{1}{\rho(\lambda,\mu)} \Big] & \lesssim \frac{2^{-p}}{2\pi} \int_0^{1-r^2} \log^p\Big( 1 + \frac{(1-r^2)^2}{(r\theta)^2}(1+o(\theta)) \Big) d\theta \\
& + \frac{2^{-p}}{2\pi} \int_{1-r^2}^{\pi}  \log^p\Big( 1 + \frac{(1-r^2)^2}{(r\theta)^2}(1+o(\theta)) \Big) d\theta \\
& = A + B.
\end{align*}
Note that in the previous estimates, for symmetry reasons we can restrict the integration to $[0,\pi]$.\\
The first integral is the one which needs more attention. Since $\log(1+x)\simeq \log(x)$ for $x$ big enough, we have for $\theta \leq 1-r^2$
$$
\log^p\Big( 1 + \frac{(1-r^2)^2}{(r\theta)^2}(1+o(\theta)) \Big) \simeq \log^p\Big(\frac{(1-r^2)^2}{(r\theta)^2}(1+o(\theta)) \Big).
$$
Furthermore note that under the same assumption 
$$
\log(1+o(\theta)) \simeq o(\theta) < 1-r^2 \simeq \log\Big(\frac{1}{r}\Big) \leq \log \Big( \frac{1-r^2}{r\theta} \Big).
$$
Using the last two estimates and the change of variable $t=(1-r^2)/(r\theta)$ we obtain that
\begin{align}
    A & \lesssim \frac{2^p}{2\pi} \int_0^{1-r^2} \log^p\Big( \frac{1-r^2}{r\theta} \Big) d\theta \nonumber \\
    & \leq \frac{2^p}{2\pi} \frac{1-r^2}{r} \int_{1}^{+\infty} \frac{\log^p(t)}{t^2} dt,
\end{align}
where the last inequality follows from $r < 1$. Since the last integral converges for every value of $p \in \bR$ and, as already mentioned, we can suppose $r\geq 1/2$, we obtain that
$$
A \leq C'_p (1-r^2). 
$$

Given that $\log(1+x) \leq x$, the second integral becomes
\begin{align*}
    B & \lesssim \frac{2^{-p}}{2\pi} \int_{1-r^2}^{\pi} \frac{(1-r^2)^2}{(r\theta)^2} d\theta \\
    & \leq \frac{2^{-p}}{2\pi} \frac{1-r^2}{r^2}\\
    & = C_p'' (1-r^2).
\end{align*}
Now define $C_p= \max (C_p',C_p'')$ and the result follows.
\end{proof}

\begin{lem} \label{off diagonal estimate}
    Let $\Lambda$ be a random sequence and $p\geq 1$. Suppose that $\lambda$ and $\mu$ belong to different dyadic annuli. Then there exists $C>0$ such that
    $$
    \bE\Big[\log^p \frac{1}{\rho(\lambda,\mu)} \Big] \leq C \min(1-|\lambda|,1-|\mu|).
    $$
\end{lem}

\begin{proof}
    Dyadic annuli have constant pseudohyperbolic width. For this reason, we can assume that if $\lambda$ and $\mu$ belong to different annuli, then there exists a constant $c>0$ such that $\rho(\lambda,\mu) \geq c$. If two points lie near the common border between adjacent annuli, we treat them as belonging to the same annulus and apply Lemma \ref{diagonal estimate}.
    Since for every $p\geq 1$ there exists $c>0$ such that $$\log^p(1/x) \leq c (1-x) \quad \text{for } \frac{1}{2} \leq x <1,$$ we have
    $$
    \log^p \frac{1}{\rho(\lambda,\mu)} = 2^{-p} \log^p \frac{1}{\rho(\lambda,\mu)^2} \lesssim 2^{-p} (1-\rho(\lambda,\mu)^2) = 2^{-p}\frac{(1-|\mu|^2)(1-|\lambda^2|)}{|1-\ol{\lambda}\mu|^2}.
    $$
    Consider $|\lambda|= r_1$ and $|\mu|= r_2$, then
    \begin{align*}
        \bE \Big[ \frac{1}{|1-\ol{\lambda}\mu|^2} \Big] & = \frac{1}{4\pi^2} \int_0^{2\pi}\int_0^{2\pi} \frac{1}{|1-r_1r_2 e^{i(t-s)}|^2} dtds \\
        & \simeq \int_0^{1} \frac{1}{(1-r_1r_2 + t)^2} dt \\
        & \lesssim  \int_0^{1-r_1r_2} \frac{1}{(1-r_1r_2)^2} dt + \int_{1-r_1r_2}^{1} \frac{1}{t^2} dt \\
        & \lesssim  \frac{1}{1-r_1 r_2}.
    \end{align*}
    Given that $1 - |\lambda||\mu| > \max (1-|\lambda|,1-|\mu|)$, we obtain
    $$
    \bE \Big[ \log^p \frac{1}{\rho(\lambda,\mu)} \Big] \lesssim \frac{(1-|\lambda|^2)(1-|\mu|^2)}{1-|\lambda||\mu|} \lesssim  \min (1-|\lambda|, 1-|\mu|).
    $$
\end{proof}

\begin{rem} \label{cochran lemma}
In \cite{cochran90}*{Lemma p.736}, the author proved that
\begin{equation} \label{cochran lem}
\bE [\log \rho(\lambda,\mu)] = \max (\log |\lambda|, \log |\mu|).
\end{equation}
Since $\log |\lambda|^{-1} \simeq (1-|\lambda|)$ for $|\lambda| \to 1$, we can rewrite \eqref{cochran lem} as
$$
\bE \Big[ \log \frac{1}{\rho(\lambda,\mu)} \Big] \simeq \min (1- |\lambda|, 1- |\mu|).
$$
From this perspective, Lemma \ref{diagonal estimate} and Lemma \ref{off diagonal estimate} together can be viewed, at least for the upper estimate, as an extension of this result, where instead of considering $\log \rho(\lambda,\mu)$, we consider $\log^p \rho(\lambda,\mu)$ for $p > 1$. We note that the proof of \eqref{cochran lem} in \cite{cochran90}*{Lemma p.736} uses Jensen’s formula to compute
$
\int \log \rho(r_1,r_2 e^{it})dt.
$
This is no longer possible  when $p>1$.

\end{rem}

We are now ready for the proof of the main result.

\begin{proof}[Proof of Theorem \ref{main: hardy-orlicz}]
We change notation for the proof of this theorem: we write $\lambda_n,\lambda_m$ instead of $\lambda,\mu$. Furthermore, without loss of generality, we can suppose that the sequence of radii is increasing, i.e. $|\lambda_n|\leq|\lambda_m|$ for $n\leq m$ and that $|\lambda_n| \geq 1/2$. Let now $\Lambda$ be a random sequence which satisfies the Blaschke condition. Consider $\lambda_n \in \Lambda$ and 
$$
B_n(z):= B_{\Lambda \setminus \{\lambda_n\}} =\prod_{k \neq n} \frac{|\lambda_k|}{\lambda_k} \frac{\lambda_k - z}{1 - \ol{\lambda_k}z}.
$$
We start proving that there exists a constant $C_p>0$ that does not depend on $n$ (but depends on $p$) such that
\begin{equation} \label{expectation log^p}
    \bE\Big[\log^p\frac{1}{B_n(\lambda_n)}\Big] \leq C_p.
\end{equation}
To this end we use Rosenthal's inequality (Theorem \ref{rosenthal ineq}) and Fatou's lemma. Consider $k\in \bN$ such that $k>p$ and $k>n$. Applying \eqref{rosenthal} with $X_m = \log \rho(\lambda_n,\lambda_m)$, $m \neq n$, we obtain
\begin{align*}
\bE\Big[ \big( \sum_{\substack{m=1\\m\neq n}}^k \log \frac{1}{\rho(\lambda_n,\lambda_m)} \big)^p \Big] & \leq 2^{p^2} \max\Big\{ \sum_{\substack{m=1\\m\neq n}}^k \bE\big[\log^p \frac{1}{\rho(\lambda_n,\lambda_m)}\big],  \big( \sum_{\substack{m=1\\m\neq n}}^k \bE\big[\log\frac{1}{\rho(\lambda_n,\lambda_m)}\big]\big)^p \Big\}.
\end{align*}
We start with the first sum. Note that from Lemma \ref{diagonal estimate} and Lemma \ref{off diagonal estimate} it follows that
$$
\bE\Big[ \log^p \frac{1}{\rho(\lambda_n,\lambda_m)} \Big] \leq C_p \min(1-|\lambda_n|, 1-|\lambda_m|).
$$
This means that
\begin{align*}
    \sum_{\substack{m=1\\m\neq n}}^k \bE \Big[ \log^p \frac{1}{\rho(\lambda_n,\lambda_m)} \Big] & \leq C_p \Big(\sum_{m=1}^{n-1} (1-|\lambda_n|) + \sum_{m=n+1}^k (1-|\lambda_m|) \Big)\\
    & = C_p \Big( (n-1)(1-|\lambda_n|) + \sum_{m=n+1}^k (1-|\lambda_m|) \Big).
\end{align*}
Taking the limit for $k \to \infty$ we obtain that
\begin{equation} \label{E [log^p]}
\sum_{\substack{m=1\\m\neq n}}^\infty \bE \Big[ \log^p \frac{1}{\rho(\lambda_n,\lambda_m)} \Big] \leq C_p \Big( (n-1)(1-|\lambda_n|) + \sum_{m=n+1}^\infty (1-|\lambda_m|) \Big).
\end{equation}
We know that $(1-|\lambda_n|)_{n \in \bN}$ is a non-increasing sequence in $\ell^1$. It is well known that this implies
$$
(n-1)(1-|\lambda_n|) \xrightarrow[n \to \infty]{} 0 \quad \text{and} \quad \sum_{m=n+1}^\infty (1-|\lambda_m|) \xrightarrow[n \to \infty]{} 0.
$$
We have actually proved that there exists a constant (that again we call) $C_p>0$ such that
$$
\sum_{\substack{m=1\\m\neq n}}^\infty \bE \Big[ \log^p \frac{1}{\rho(\lambda_n,\lambda_m)} \Big] \leq C_p.
$$
Since the above proof works also for $p=1$, we have $\bE[ \log |B_n(\lambda_n)|^{-1}] \leq C_1$ and then
$$
\Big( \sum_{\substack{m=1\\m\neq n}}^{\infty} \bE\Big[\log\frac{1}{\rho(\lambda_n,\lambda_m)}\Big]\Big)^p \leq C_1^p.
$$
This gives that
$$
\bE\Big[ \big( \sum_{\substack{m=1\\m\neq n}}^k \log \frac{1}{\rho(\lambda_n,\lambda_m)} \big)^p \Big] \leq 2^{p^2} \max(C_p, C_1^p).
$$
By Fatou’s lemma (and observing that the partial sums are increasing so that limits are inferior limits) we get
\begin{align*}
    \bE \Big[ \log^p \frac{1}{B_n(\lambda_n)} \Big] & = \bE\Big[ \big( \lim_{k \to \infty} \sum_{\substack{m=1\\m\neq n}}^k \log \frac{1}{\rho(\lambda_n,\lambda_m)} \big)^p \Big] \\
    & \leq \lim_{k \to \infty} \bE\Big[ \big( \sum_{\substack{m=1\\m\neq n}}^k \log \frac{1}{\rho(\lambda_n,\lambda_m)} \big)^p \Big] \\
    & \leq 2^{p^2} C_p.
\end{align*}

Define $X= \sum_{n \in \bN} (1-|\lambda_n|) \log^p|B_n(\lambda_n)|^{-1}$. Note that $X$ is a sum of non-negative random variables. If we prove that its expectation is finite, then $X < \infty$ almost surely. We have that
\begin{align*}
    \bE [X] & = \bE \Big[ \sum_{n \in \bN} (1-|\lambda_n|) \log^p\frac{1}{|B_n(\lambda_n)|} \Big] \\
    & = \sum_{n \in \bN} (1-|\lambda_n|) \bE\Big[ \log^p\frac{1}{|B_n(\lambda_n)|}\Big] \\
    & \leq 2^{p^2}C_p \sum_{n \in \bN} (1-|\lambda_n|) < \infty.
\end{align*}
This implies that $\bP(X < \infty)=1$. The statement now follows from Proposition \ref{generalization H_p}. 

Since every free interpolating sequence is a zero set and zero sets for $\cH_p$ are Blaschke sequences, the necessary part follows.
\end{proof}

\begin{rem}
    Carleson sequences are in general not free interpolating for the Nevanlinna class, neither for the Smirnov class or for $\cH_p$. Indeed, a Carleson sequence is a finite union of $H^{\infty}$-interpolating and thus weakly separated sequences. Now, it is not difficult to construct a union of two sequences coming close together in the pseudohyperbolic metric  and that cannot be interpolating for $\cN$ (and thus neither for $\cN^+$ nor for $\cH_p$).  Pick, for instance, $\Lambda_1=(1-2^{-n})_{n \in \bN}$ (which is interpolating for $H^{\infty}$) and $\Lambda_2=(\mu_n)_{n \in \bN}$ with $\rho(1-2^{-n},\mu_n)=e^{-n2^n}$ (so that $\Lambda_2$ is also interpolating for $H^{\infty}$). Then, defining $\Lambda=\Lambda_1\cup\Lambda_2$ (which is a Carleson sequence) and $B=B_\Lambda$, we have
    $$
    \log1/|B_n(\lambda)|\geq \log 1/\rho(1-2^{-n},\mu_n)=n2^n > P[\mu](1-2^{-n})
    $$ 
    for every positive finite measure $\mu$ on $\bT$. Indeed, given any positive harmonic function $h$, we have $h(1-2^{-n})\le C 2^n$.
    
    The situation changes drastically when we consider random sequences. In fact every almost surely Carleson sequence is necessary a Blaschke sequence, which means that it is almost surely free interpolating for the Nevanlinna or Smirnov classes.
\end{rem}

We conclude this paper observing that the random model we consider is not adapted to the interpolation problem solved by Naftalevi\v{c}.

\begin{prop} \label{naftalevic:random}
    Let $\Lambda$ be a random sequence. Then
    $$
    \bP \big(\Lambda \text{ is interpolating for $\cN$ in the sense of Naftalevi\v{c}} \big)=0.
    $$
\end{prop}

\begin{proof}
As shown in \cite{cochran90}, all the elements of $\bT$ are almost surely accumulation points for the random sequence $\Lambda$. In particular it follows that with probability $0$ the whole sequence is contained in a finite union of Stolz angles. Since this is a necessary condition for $\Lambda$ to be interpolating in the sense of Naftalevi\v{c}, the result follows.
\end{proof}

\section*{Acknowledgements}
I would like to thank my supervisor, professor Andreas Hartmann, for introducing me to the problem and for the many fruitful discussions we had. In addition, I would like to thank Nikolaos Chalmoukis for some useful comments on Lemma \ref{diagonal estimate}.

\bibliographystyle{habbrv}
\bibliography{biblio}

%\printbibliography

\end{document}